\documentstyle[12pt]{amsart}
\newtheorem{theorem}{Theorem}[section]

\newtheorem{proposition}[theorem]{Proposition}

\numberwithin{equation}{section}

\begin{document}

\title [Symplectic and Lagrangian mean curvature flows]
 {The second type Singularity of symplectic and Lagrangian mean curvature flows}
\author{Xiaoli Han, Jiayu Li}

\address{Math. Group, The abdus salam ICTP\\ Trieste 34100,
   Italy}
\email{xhan@@ictp.it}

\address{Math. Group, The abdus salam ICTP\\ Trieste 34100,
   Italy\\
   and Academy of Mathematics and Systems Sciences\\ Chinese Academy of
Sciences\\ Beijing 100080, P. R. of China. } \email{jyli@@ictp.it}


\date{}

\maketitle

\begin{abstract}
In this paper we mainly study the type II singularities of the
mean curvature flow from a symplectic surface or from an almost
calibrated Lagrangian surface in a K\"ahler surface. We study the
relation between the maximum of the K\"ahler angle and the maximum
of $|H|^2$ on the limit flow.
\end{abstract}


\section{Introduction}
In this paper, we continue to study the symplectic mean curvature
flow and Lagrangian mean curvature flow (\cite{CL1}, \cite{CL2},
\cite{CL3} \cite{HL1}, \cite{HL2}, \cite{Sm1}, \cite{Wa}) in a
K\"ahler surface. Suppose $M$ is a compact K\"ahler surface. Let
$\Sigma$ be a smooth surface in $M$ and $\omega$, $\langle\cdot,
\cdot\rangle$ be the K\"ahler form and the K\"ahler metric on $M$
respectively. The K\"ahler angle $\alpha$ of $\Sigma$ in $M$ is
defined by \cite{CW}
$$ \omega|_\Sigma=\cos\alpha d\mu_\Sigma
$$ where $d\mu_\Sigma$ is the area element of $\Sigma$ of the
induced metric from $\langle, \rangle$. We call $\Sigma$ a {\it
symplectic} surface if $\cos\alpha>0$, a {\it Lagrangian} surface
if $\cos\alpha=0$, a {\it holomorphic curve} if $\cos\alpha=1$. In
addition, we assume that $M$ is a Calabi-Yau manifold of complex
dimension $2$ with a complex structure $J$, i.e, a K3 surface. We
consider a parallel holomorphic $(2, 0)$ form,
$$ \Omega=dz_1\wedge dz_2. $$
If a surface $\Sigma$ is Lagrangian then (see \cite{HaL})
$$\Omega |_\Sigma=e^{i\theta}d\mu_\Sigma,
$$
where $\theta$ is a  multivalued function called Lagrangian angle.
If $ \cos\theta> 0, $ then $\Sigma$ is called {\it almost
calibrated}. If $\theta=costant$, then $\Sigma$ is called {\it
special} Lagrangian.

It is proved in \cite{CL1} and \cite{Wa} that, if the initial
surface is symplectic, then along the mean curvature flow, at each
time $t$ the surface is still symplectic. Thus we speak of
symplectic mean curvature flow. It is proved in \cite{Sm1},
\cite{Sm2} that, if the initial surface is Lagrangian, then along
the mean curvature flow, at each time $t$ the surface is still
Lagrangian. Thus we speak of Lagrangian mean curvature flow.

In \cite{HL1} we showed that, if the scalar curvature of the
compact K\"ahler-Einstein surface $M$ is positive and the initial
surface is sufficiently close to a holomorphic curve, then the
mean curvature flow has a global solution and converges to a
holomorphic curve.

In general, the mean curvature flow may produce singularities. The
beautiful results on the nature of singularities of the mean
curvature flow of convex hypersurfaces have been obtained by
Huisken-Sinestrari \cite{HS1}, \cite{HS2} and White \cite{W}. For
symplectic mean curvature flow, Chen-Li \cite{CL1} and Wang
\cite{Wa} proved that there is no Type I singularity. At a Type II
singular point, Chen-Li \cite{CL2}, \cite{CL3} proved that, the
rescaled surfaces converge weakly (in the sense of measure) to a
stationary tangent cone which is flat.

If we consider the strong convergence of the rescaled surfaces
$\Sigma^k_s$ in $B_R(0)$ around a type II singular point, let
$|A_k|$ be the second fundamental forms of $\Sigma^k_s$ in
$B_R(0)$, then we have that $|A_k|^2\leq 4$ in $B_R(0)$ during the
rescaling process. Thus by Arzela-Ascoli theorem,
$\Sigma^k_s\to\Sigma^\infty_s$ in $C^2(B_R(0)\times [-R, R])$ for
any $R>0$ and any $B_R(0)\subset {\mathbb {C}}^2$. By the
definition of the type II singularity, we know that
$\Sigma^\infty_s$ is defined on $(-\infty, +\infty)$ and
$\Sigma_s^\infty$ also evolves along the mean curvature flow in
${\mathbb{C}}^2$ with the Euclidean metric. We call
$\Sigma^\infty_s$ the limit flow at $X_0$. See Section 2 for
details.

In this paper, we mainly study the nature of the limit flow
$\Sigma^\infty_s$. For this purpose, we consider a general mean
curvature flow $\Sigma_t$ in ${\mathbb{R}}^4$ which exists
globally with bounded second fundamental forms. In particular,
translating soliton to the mean curvature flow is a special case.
Recently in \cite{HL2} we proved that there is no translating
soliton with $\cos\alpha\geq\delta$ to the symplectic mean
curvature flow or to the almost calibrated Lagrangian mean
curvature flow where $\delta>0$ is a constant depending only on
the speed of the soliton. Since $\Sigma_t$ come from the blow up,
it is natural to assume that on $\Sigma_t$, we have
\begin{equation}\label{e0.0} cR^2\leq\mu_t(\Sigma_t\cap B_R(0))\leq
CR^2,
\end{equation} where $0<c<C<\infty$ are constants which are
independent of $t$ and $R$.

\vspace{.2in}

\noindent {\bf Main Theorem 1} {\it Suppose that $\Sigma_t$, $t\in
(-\infty, 0]$ is a complete symplectic mean curvature flow in
${\mathbb{C}}^2$ which satisfies (\ref{e0.0}). Assume that
$\sup_{t\in (-\infty, 0]}\sup_{\Sigma_t}|A|^2= 1$. If
$h^2=\sup_{t\in (-\infty, 0]}\sup_{\Sigma_t} |H|^2$ and
$\delta=\inf_{t\in (-\infty, 0]}\inf_{\Sigma_t}\cos\alpha$, then
$\delta e^{\frac{h^2}{4}}\leq 1$.}

\vspace{.1in}

Analogously in the almost calibrated Lagrangian mean curvature
flow, we have \vspace{.2in}

\noindent {\bf Main Theorem 2} {\it Suppose that $\Sigma_t$, $t\in
(-\infty, 0]$ is a complete almost calibrated Lagrangian mean
curvature flow in ${\mathbb{C}}^2$ which satisfies (\ref{e0.0}).
Assume further that $\sup_{t\in (-\infty, 0]}\sup_{\Sigma_t}|A|^2=
1$. If $h^2=\sup_{t\in (-\infty, 0]}\sup_{\Sigma_t} |H|^2$ and
$\delta=\inf_{t\in (-\infty, 0]}\inf_{\Sigma_t} \cos\theta$, then
$\delta e^{\frac{h^2}{2}}\leq 1$.} \vspace {.1in}

The authors would like to thank the referees for their valuable
comments which improved this paper very much.
\section{Preparations}

In this section we define the rescaled surfaces and study the
strong convergence of the rescaled sequence at a type II singular
point, which is more or less standard. However we can not find it
in a reference, so we give all details here. It may be interesting
in its own right. Suppose that $T$ is discrete singular time, that
means there exists $\varepsilon>0$ such that the mean curvature
flow is smooth in $[T-\varepsilon, T)$. Assume that $(X_0, T)$ is
a type II singular point of the mean curvature flow in $M$. Since
this is type II singularity, then for any sequence $\{r_k\}$ with
$r_k\to 0$,
\begin{eqnarray*}&&\max_{\sigma\in (0, r_k/2]}\sigma^2 \max_{[T-(r_k-\sigma)^2,
T-(r_k/2)^2]} \max_{\Sigma_t\cap B_{r_k-\sigma}(X_0)}|A|^2\\
&\geq& (r_k/2)^2\max_{\Sigma_{T-(r_k/2)^2}\cap B_{r_k/2}(X_0)}|A|^2\\
&=&(T-(T-(r_k/2)^2))\max_{\Sigma_{T-(r_k/2)^2}\cap
B_{r_k/2}(X_0)}|A|^2\\ &\rightarrow& +\infty
\end{eqnarray*}
We choose $\sigma_k\in (0, r_k/2]$ such that
$$\sigma_k^2 \max_{[T-(r_k-\sigma_k)^2, T-(r_k/2)^2]}
\max_{\Sigma_t\cap B_{r_k-\sigma_k}(X_0)}|A|^2=\max_{\sigma\in (0,
r_k/2]}\sigma^2 \max_{[T-(r_k-\sigma)^2, T-(r_k/2)^2]}
\max_{\Sigma_t\cap B_{r_k-\sigma}(X_0)}|A|^2.
$$
Let $t_k\in [T-(r_k-\sigma_k)^2, T-(r_k/2)^2]$ and
$F(x_k,t_k)=X_k\in\bar{B}_{r_k-\sigma_k}(X_0)$ satisfy
$$\lambda_k^2=|A|^2(X_k)=|A|^2(x_k, t_k)=\max_{[T-(r_k-\sigma_k)^2, T-(r_k/2)^2]}
\max_{\Sigma_t\cap B_{r_k-\sigma_k}(X_0)}|A|^2.$$ Obviously, we
have $(X_k, t_k)\to (X_0, T)$ and
$\lambda_k^2\sigma_k^2\to\infty$. In particular,
\begin{equation}\label{e0.2} \max_{[T-(r_k-\sigma_k/2)^2, T-(r_k/2)^2]} \max_{\Sigma_t\cap
B_{r_k-\sigma_k/2}(X_0)}|A|^2\leq 4\lambda_k^2, \end{equation} and
hence \begin{equation}\label{e0.3}\max_{[t_k-(\sigma_k/2)^2, t_k]}
\max_{\Sigma_t\cap B_{r_k-\sigma_k/2}(X_0)}|A|^2\leq 4\lambda_k^2.
\end{equation}

We now describe the rescaling process around $(X_0,T)$ in details.
The argument is discussed with J. Chen. In the following we denote
the points of the image of $F$ or $F_k$ in $M$ by capital letters.
We choose a normal coordinates in $B_r(X_0)$ using the exponential
map, where $B_r(X_0)$ is a metric ball in $M$ centered at $X_0$
with radius $r$ ($0<r<i_M/2$, $i_M$ is the injective radius of
$M$). We express $F$ in its coordinates functions. Consider the
following sequences,
\begin{equation}\label{e0.6} F_k(x, s)=\lambda_k(F(x_k+x,
t_k+\lambda_k^{-2}s)-F(x_k, t_k)),~~~~~~ s\in
[-\lambda_k^2\sigma_k^2/4, \lambda_k^2(T-t_k)].\end{equation} We
denote the rescaled surfaces by $\Sigma^k_s$ in which $d\mu^k_s$
is the induced area element from $M$. For any $R>0$, let $B_R(0)$
be a ball in ${\mathbb{R}}^4$ with radius $R$ in the Euclidean
metric and centered at $0$. Then
$$\Sigma^k_s\cap  B_R(0)=\{|F_k(x, s)|\leq R\}, $$ it is clear
that for any fixed $R>0$, $\lambda_k^{-1}R<r/2, r_k<r/2$ as $k$
sufficiently large, then the surface $\Sigma^k_s$ is defined on
$B_R(0)$ because
\begin{eqnarray*} \exp_{X_0} (\lambda_k^{-1}\{|F_k(x, s)|\leq
R\})&\subset& \exp_{X_0} (|F-X_0|\leq \lambda_k^{-1}R+r_k)\\
&\subset& B_{\lambda_k^{-1}R+r_k}(X_0)\subset B_{r}(X_0).
\end{eqnarray*}  Moreover, we pull back the metric on $B_r(X_0)\subset
M$ via $\exp_{X_0}$ so that we get a metric $h$ on the Euclidean
ball $B_r(0)$. Then for any fixed $R>0$ such that
$\lambda_k^{-1}R<r/2$, we can define a metric $h_{k, R}$ on
$B_R(0)$,
$$( h_{k, R})_{ij}(X)=\lambda_k^{2}h(\lambda_k^{-1}X+X_k).
$$ With respect to this metric $\Sigma^k_s$ evolves along the mean
curvature flow, which will be derived as follows.

If $g^k_s$ is the metric on $\Sigma^k_s$ which is induced from the
metric $g(\cdot, t_k+\lambda^{-1}_k s)$ on
$\Sigma_{t_k+\lambda_k^{-1}s}$, it is clear that
$$ (g^k_s)_{ij}(X)=\lambda_k^2 g_{ij}(\lambda_k^{-1}X+X_k,
t_k+\lambda_k^{-2}s),
$$ and
$$ (g^k_s)^{ij}(X)=\lambda_k^{-2} g^{ij}(\lambda_k^{-1}X+X_k,
t_k+\lambda_k^{-2}s).
$$ In this setting $(\Sigma^k_s, g_s^k)$ is an isometric immersion in
$(B_R(0), h_{k, R})$. Let $A_k$, $H_k$ be the second fundamental
form and the mean curvature vector of $(\Sigma^k_s, g^k_s)$ in
$(B_R(0), h_{k, R})$ respectively. Let $\bar\Gamma^k$,
$\Gamma^k_s$ be the Christoffel symbols of $h_{k, R}$ on $B_R(0)$
and the Christoffel symbols of $g^k_s$ on $\Sigma^k_s$ . Since
$F_k$ is an isometric immersion in $(B_R(0), h_{k, R})$ with
respect to the induced metric, hence by the Gauss equation we
have,
\begin{eqnarray}\label{e0.4}
(A_k)_{ij} &=& \sum_{\alpha=1, 2}(h_k)_{ij}^\alpha \nu^k_{s\alpha}
\nonumber\\ &=& -\partial^2_{ij} F_k +\sum_{l=1,
2}(\Gamma^k_s)^l_{ij}\partial_l F_k-\sum_{\alpha,\beta,\gamma=1,
4}(\bar\Gamma^k)^\alpha_{\beta\gamma}\partial_i F^\beta_k
\partial_j F^\gamma_k \nu^k_{s\alpha},
\end{eqnarray} where $\{\nu^k_{s\alpha},\alpha=1, 2\}$ are bases of the normal
space of $\Sigma^k_s$ in $(B_R(0), h_{k, R})$. Let
$\Gamma_{t_k+\lambda_k^{-2}s}$ be the Christoffel symbols on
$\Sigma_{t_k+\lambda_k^{-2}s}$ and $\bar\Gamma$ be the Christoffel
symbols on $M$. It is not hard to check that
$$\bar\Gamma^k(X)=\bar\Gamma(\lambda_k^{-1}X+X_k),~~~~~~~~~~~
\Gamma^k_s(X)=\Gamma_{t_k+\lambda_k^{-2}s}(\lambda_k^{-1}X+X_k).
$$ Thus from (\ref{e0.4}), we get that,
\begin{eqnarray}\label{e0.5}
(A_k)_{ij} &=& \lambda_k(-\partial^2_{ij} F+\sum_{l=1,
2}(\Gamma_{t_k+\lambda_k^{-2}s})^l_{ij}\partial_l
F_k-\sum_{\alpha,\beta,\gamma=1,
4}\bar\Gamma^\alpha_{\beta\gamma}\partial_i F^\beta_k
\partial_j F^\gamma_k \nu_\alpha )\nonumber\\ &=& \lambda_k A_{ij},
\end{eqnarray} where $\{v_\alpha, \alpha=1, 2\}$ are bases of the normal space of
$\Sigma_{t_k+\lambda_k^{-2}s}$ in $M$.  Therefore,
\begin{eqnarray*}
|A_k|^2 &=&\lambda_k^{-2}|A|^2,\\
H_k&=&\lambda_k^{-1}H,\\ |H_k|^2 &=&\lambda_k^{-2}|H|^2.
\end{eqnarray*}

Set $t=t_k+\lambda_k^{-2}s$, it is easy to check that
\begin{eqnarray*}
\frac{\partial F_k}{\partial s}&=&\lambda_k^{-1}\frac{\partial
F}{\partial t}.
\end{eqnarray*}

Therefore, it follows that the scaled surface also evolves by a
mean curvature flow
\begin{equation}\label{e0.1} \frac{\partial F_k}{\partial s}=H_k
\end{equation} in $B_{\lambda_k\sigma_k}(0)$, where $s\in [-\lambda_k^2\sigma_k^2/4,
\lambda_k^2(T-t)].
$

By (\ref{e0.2}) and (\ref{e0.3}) we know that,
$$|A_k|(0, 0)=1,~~~~~~~~~|A_k|^2\leq 4$$ in $B_{\lambda_k\sigma_k}(0)$
and $s\in [-\lambda_k^2\sigma_k^2/4, \lambda_k^2(T-t)]. $ Since
$(X_0, T)$ is a type II singularity, then
$\lambda_k^2\sigma_k^2\to\infty$ and
$\lambda_k^2(T-t_k)\to\infty$.  Thus by Arzela-Ascoli theorem,
$\Sigma^k_s\to\Sigma^\infty_s$ in $C^2(B_R(0)\times [-R, R])$ for
any $R>0$ and any $B_R(0)\subset {\mathbb {C}}^2$. By
(\ref{e0.6}), we know that $\Sigma^\infty_s$ is defined on
$(-\infty, +\infty)$. Since for each fixed $R>0$,
$\lambda_k^{-1}X+X_k\to X_0$ for $X\in B_R(0)$ as $k\to \infty$,
then $h_{k, R}$ converges uniformly in $B_R(0)$ to the Euclidean
metric as $k\to \infty$, and the Christoffel symbols
$(\bar\Gamma^k)$ of $h_{k, R}$ converges uniformly in $B_R(0)$ to
$0$ as $k\to\infty$, we see that $\Sigma^\infty_s$ also evolves
along the mean curvature flow in ${\mathbb{C}}^2$ with the
Euclidean metric. We call $\Sigma^\infty_s$ the limit flow at
$X_0$.

In the rest part of this section, we estimate the different of
$A_k, H_k$ and $A_k^0, H_k^0$ where $A^0_k$ and $H^0_k$ are the
second fundamental form and the mean curvature vector of
$\Sigma^k_s$ in the Euclidean metric on $B_R(0)$ respectively.
Although it is not needed in this paper, it is interesting in its
own right.

 Let $\Gamma^{0k}_s$
be the Christoffel symbols of $\Sigma^k_s$ for the Euclidean
metric on $B_R(0)$ and $\{\nu^{0k}_{s\alpha}: \alpha=1, 2\}$ be
bases of the normal space of $\Sigma^k_s$ with respect to the
Euclidean metric on $B_R(0)$ . Similarly, considering $F_k$ as an
isometric immersion in $B_R(0)$ with Euclidean metric, we have,
\begin{eqnarray}\label{e0.5}
(A^0_k)_{ij} &=& \sum_{\alpha=1,
2}(h_0)^\alpha_{ij}(\nu^{0k}_s)_\alpha=-\partial^2_{ij}
F_k+\sum_{l=1, 2}(\Gamma^{0k}_s)^l_{ij}\partial_l F_k.
\end{eqnarray}

Note that the induced metric on $\Sigma^k_s$ from $h_{k, R}$ is
given by $\langle\partial F_k, \partial F_k\rangle_{h_{k, R}}$, so
it holds
$$|\partial F_k|^2_{h_{k, R}}=2,
$$
which in turn implies that for $k$ sufficiently large and $R$
fixed $|\partial F^\alpha_k|$ is uniformly bounded in $B_R(0)$
with Euclidean metric.

Using the Euclidean metric on $B_R(0)$, we decompose the tangent
bundle of $B_R(0)$ along $\Sigma^k_s$ into the tangential
component $T\Sigma^k_s$ and the normal component
$T^\perp\Sigma^k_s$. Let $A^\perp_k: T\Sigma^k_s\times
T\Sigma^k_s\to T^\perp\Sigma^k_s$ be the normal component of
$A_k$. Notice that $A_k^\perp-A^0_k$ lies in $T^\perp\Sigma^k_s$
and $\partial_i F_k$ lies in $T\Sigma^k_s$, it follows from
(\ref{e0.4}) and (\ref{e0.5}) that,
$$\sup_{B_R(0)}|A_k^\perp-A_k^0|\leq
C\sup_{B_R(0)}|\bar\Gamma^k|\to 0
$$ as $k\to\infty$ for any fixed $R>0$. From the uniform
convergence of the metrics $h_{k, R}$ to the Euclidean metric,
$$|A_k^\perp|\leq |A_k|\leq 2|A_k|_{h_{k, R}}
$$ for any fixed $R>0$ and sufficiently large $k$. Hence, there
exist positive constants $\delta_{k, R}$ which tend to $0$ as
$k\to \infty$ such that
$$|A^0_k|=|A_k^\perp|+\delta_{k, R}\leq 2|A_k|_{h_{k,
R}}+\delta_{k, R}
$$ for all sufficiently large $k$ and any fixed $R>0$; and
similarly there exist constants $\delta'_{k, R}>0$ with
$\delta'_{k, R}\to 0$ as $k\to\infty $ such that
$$ |H^0_k|\leq 2|H_k|_{h_{k, R}}+\delta'_{k, R}
$$ for sufficiently large $k$ and for any given $R>0$.

\section{Proof of the Main Theorems}

Now we begin to prove our Main Theorems. We first prove Main
Theorem $2$. Let $H(X, X_0, t, t_0)$ be the backward heat kernel
on ${\mathbb{R}}^4$. Let $\Sigma_t$ be a smooth family of surfaces
in ${\mathbb{R}}^4$ defined by $F_t: \Sigma\to {\mathbb{R}}^4$.
Define
$$\rho(X,  t)=(4\pi (t_0-t))H(X, X_0, t, t_0)=\frac{1}{ 4\pi
(t_0-t)} \exp{-\frac{|X-X_0|^2}{4(t_0-t)}}
$$ for $t<t_0$, such that
$$\frac{d}{dt}\rho=-\Delta\rho-\rho\left(\left|H+\frac{(X-X_0)^\perp}{2(t_0-t)}
\right|^2-|H|^2\right).
$$ where $(X-X_0)^\perp$ is the normal component of $X-X_0$.

 Define
$$\Psi_{X_0, t_0}(X, t)=\int_{\Sigma_t}\frac{1}{\cos\theta}\rho(X, t)d\mu_t.
$$

\allowdisplaybreaks
\begin{proposition}\label{mono1}
Along the almost calibrated Lagrangian mean curvature flow
$\Sigma_t$ in ${\mathbb{C}}^2$, we have,
\begin{eqnarray*}
&&\frac{\partial}{\partial t}\Psi_{X_0, t_0}(X, t)\nonumber\\&&=
-\left( \int_{\Sigma_t}\frac{1}{\cos\theta}\rho (F, t)
\left|H+\frac{(F-X_0)^{\perp}}{2(t_0-t)}\right|^2d\mu_t \right.\nonumber \\
&&\left.+ \int_{\Sigma_t}\frac{1}{\cos\theta}\rho (F,
t)|H|^2d\mu_t +\int_{\Sigma_t}\frac{2}{\cos^3\theta}\left|\nabla
\cos\theta\right|^2\rho (F, t)d\mu_t\right).
\end{eqnarray*}
\end{proposition}

{\it Proof.} From  the evolution equation of Lagrangian angle
(\cite{Sm1}, \cite{Sm2}),
\begin{equation}\label{e3}(\frac{\partial}{\partial
t}-\Delta)\cos\theta=|H|^2\cos\theta,\end{equation} we know that
\begin{eqnarray}\label{e9}
(\frac{\partial}{\partial
t}-\Delta)\frac{1}{\cos\theta}=-\frac{|H|^2}{\cos^2\theta}
-2\frac{|\nabla\cos\theta|^2}{\cos^3\theta}.
\end{eqnarray} Recall the general formula $(7)$ in
\cite{EH}, for a smooth function $f=f(x, t)$ on $\Sigma_t$ with
polynomial growth at infinity,
\begin{equation}\label{e8}
\frac{d}{dt}\int_{\Sigma_t}f\rho d\mu_t= \int_{\Sigma_t}
(\frac{d}{dt}f-\Delta f)\rho d\mu_t-\int_{\Sigma_t}f\rho
\left|H+\frac{(X-X_0)^\perp}{2(t_0-t)}\right|d\mu_t.
\end{equation}

Choosing $f=\frac{1}{\cos\theta}$ in (\ref{e8}) and putting
(\ref{e9}) into (\ref{e8}), we get our monotonicity formula.

\hfill Q. E. D.

\vspace{.1in}

{\it Proof of Main Theorem 2.} If $h=0$, or $\delta=0$, or
$\delta=1$, it is evident that the result holds. Now we assume
that $h>0$, $0<\delta<1$ and argue it by contradiction. Suppose
that $\delta> e^{-\frac{h^2}{2}}$. Fix $R>0$. First we claim that
there exists a sequence $\{s_i\}$ such that $s_i\to -\infty$ as
$i\to \infty$ and $\lim_{i\to \infty}\max_{\Sigma_{s_i}\cap
B_R(X_0)}|H|^2=0$. Without loss of generality, we assume $X_0=0$.
Integrating the monotonicity formula in Proposition 2.1 with
$t_0=0$ from $2s$ to $s$ for $s<0$, we get that,
\begin{eqnarray*}
&&\int_{\Sigma_{2s}}\frac{1}{\cos\theta(x,
2s)}\frac{1}{-2s}e^{\frac{|F^\perp|^2}{2s}}d\mu_{2s}-\int_{\Sigma_s}\frac{1}{\cos\theta(x,
s)}\frac{1}{-s}e^{\frac{|F^\perp|^2}{s}}d\mu_s\\ &&\geq
\int_{2s}^s\int_{\Sigma_t}\frac{1}{\cos\theta}\rho(F,
t)|H|^2d\mu_t dt.
\end{eqnarray*}  By Proposition \ref{mono1},
we know that $\int_{\Sigma_s}\frac{1}{\cos\theta}\rho(F, s)$ is
nonincreasing as $s$. Since $\cos\theta$ is bounded below by
$\delta$, for any $t<0$,
\begin{eqnarray*}
\int_{\Sigma_t}\frac{1}{\cos\theta}\rho(X, t) d\mu_t  &\leq&
1/\delta\int_{\Sigma_t} \rho(X, t) d\mu_t \\ &\leq& C/\delta
\int_0^\infty\int_{\Sigma_t\cap \partial
B_\rho(0)}\frac{1}{0-t}e^{\frac{\rho^2}{t}}d\sigma_t d\rho \\
&\leq& \frac{C}{-t}\int_0^\infty
e^{\frac{\rho^2}{t}}\frac{d}{d\rho}vol(B_\rho(0)\cap
\Sigma_t)d\rho \\
&\leq&\frac{C}{-t}[e^{\frac{\rho^2}{t}}vol(\left. B_\rho(0)\cap
\Sigma_t)\right|_{\rho=0}^\infty-\int_0^\infty
vol(B_\rho(0)\cap\Sigma_t)e^{\frac{\rho^2}{t}}\frac{2\rho}{t}d\rho],
\end{eqnarray*} where we denote by $C>0$ the constants which does
not depend on $t$ and may change from one line to another line.
Since we have assumed that $cR^2\leq B_R(0)\cap \Sigma_t\leq CR^2$
in (\ref{e0.0}),  thus we have,
\begin{eqnarray*}
\int_{\Sigma_t}\frac{1}{\cos\theta}\rho(X, t) d\mu_t &\leq& C
[\left.\frac{1}{-t}e^{\frac{\rho^2}{t}}\rho^2\right|_{\rho=0}^\infty+\int_0^\infty
\frac{2\rho^3}{t^2}e^{\frac{\rho^2}{t}}d\rho]\\ &\leq& C
[\left.\frac{1}{-t}e^{\frac{\rho^2}{t}}\rho^2+e^{\frac{\rho^2}{t}}\frac{\rho^2}{t}
-e^{\frac{\rho^2}{t}}]\right|_{\rho=0}^\infty \\ &\leq & C.
\end{eqnarray*}  Thus the quantity
$\int_{\Sigma_s}\frac{1}{\cos\theta}\rho(F, s)$ is uniformly
bounded above. Moreover, by the mean value theorem there is
$s'\in [2s, s]$ such that,
\begin{eqnarray*}
&&\int^{s}_{2s}\int_{\Sigma_t}\frac{1}{\cos\theta}\frac{1}{-t}
e^{\frac{|F|^2}{t}}|H|^2 d\mu_t\\
&&=-s\int_{\Sigma_{s'}}\frac{1}{\cos\theta}\frac{1}{-s'}
e^{\frac{|F|^2}{s'}}|H|^2 d\mu_{s'}\\
&&\geq C e^{\frac{R^2}{s}}\int_{\Sigma_{s'}\cap B_R(0)} |H|^2
d\mu_{s'},
\end{eqnarray*}  where $C$ is independent of $s$. Thus we can find
a sequence $\{s_i\}$ such that $s_i\to-\infty$ as $i\to\infty$ and
$$\int_{\Sigma_{s_i}\cap B_R(0)} |H|^2 d\mu_{s_i}\to 0~~~~~~~{\rm
as}~~~~~~~~ i\to\infty.
$$ Since the second fundamental forms of $\Sigma_{s_i}$ are bounded above and
$\Sigma_{s}$ satisfy the mean curvature flow equation, then
$\Sigma_{s_i}$ strongly converges to a smooth limit surface
$\Sigma_{-\infty}$ in $B_R(0)$. Therefore,
\begin{eqnarray}\label{e7}\lim_{i\to \infty}\max_{\Sigma_{s_i}\cap B_R(0)}|H|^2=0.
\end{eqnarray} The identity can also be proved by Morse iteration.

Now we use gradient estimate to prove our theorem. For this
purpose we introduce a new function $f(X,
t)=\frac{e^{p|H|^2}}{\cos^2\theta}$, where $t\in [s_i, 0]$,
$\{s_i\}$ is the sequence in (\ref{e7}), and $p$ is constant such
that $1-p>0$. \allowdisplaybreaks
\begin{eqnarray*}\label{e4}
(\Delta-\frac{\partial}{\partial
t})f&=&\frac{1}{\cos^2\theta}(\Delta-\frac{\partial}{\partial
t})e^{p|H|^2}+e^{p|H|^2}(\Delta-\frac{\partial }{\partial
t})\frac{1}{\cos^2\theta}\nonumber\\&&+2\nabla
e^{p|H|^2}\cdot\nabla\frac{1}{\cos^2\theta}.
\end{eqnarray*} Using the evolution equation for $|H|^2$ in
${\mathbb{R}}^4$:
$$(\Delta-\frac{\partial}{\partial t}) |H|^2=2|\nabla H|^2-2(H^\alpha
h^\alpha_{ij})^2, $$ we get
\begin{eqnarray*}
(\Delta-\frac{\partial}{\partial t}) e^{p |H|^2}
&=&e^{p|H|^2}(4p^2|H|^2|\nabla |H||^2+2p|\nabla H|^2-2p|H^\alpha
h^\alpha_{ij}|^2)\\ &\geq& e^{p|H|^2}(4p^2|H|^2|\nabla
|H||^2+2p|\nabla H|^2-2p|H|^2|A|^2)\\ &\geq&
e^{p|H|^2}(4p^2|H|^2|\nabla |H||^2+2p|\nabla H|^2-2p|H|^2).
\end{eqnarray*} Since
\begin{eqnarray*}
\nabla e^{p |H|^2}  &=&\nabla (f\cos^2\theta) \\ &=&
\cos^2\theta\nabla f+2f\cos\theta\nabla\cos\theta,
\end{eqnarray*} we have,
\begin{eqnarray*}
\nabla e^{p |H|^2}\cdot\nabla\frac{1}{\cos^2\theta}
&=&\cos^2\theta\nabla
f\cdot\nabla\frac{1}{\cos^2\theta}-\frac{4f}{\cos^2\theta}|\nabla\cos\theta|^2.
\end{eqnarray*} Using the evolution equation (\ref{e3}) we get,
$$(\Delta-\frac{\partial}{\partial t})\frac{1}{\cos^2\theta}=\frac{6|\nabla\cos\theta|^2}
{\cos^4\theta}+\frac{2|H|^2}{\cos^2\theta}.
$$ So,
\begin{eqnarray}\label{e10}
(\Delta-\frac{\partial}{\partial t})f &\geq& f(4p^2|H|^2|\nabla
|H|^2|+2p|\nabla
H|^2+2(1-p)|H|^2-2\frac{|\nabla\cos\theta|^2}{\cos^2\theta})
\nonumber\\ &&+2\cos^2\theta\nabla
f\cdot\nabla\frac{1}{\cos^2\theta}.
\end{eqnarray}

Let $\psi(r)$ be a $C^2$ function on $[0, \infty)$ such that
\begin{eqnarray*}
\psi(r)=\left\{\begin{array}{clcr}  1 &{\rm if} & r\in [0, \frac{1}{2}]\\
0 &{\rm if}  & r\geq 1
\end{array}\right.
\end{eqnarray*}
$$ 0\leq\psi(r)\leq 1, \psi'(r)\leq 0, \psi''(r)\geq -C ~~~~~{\rm and
}~~~~~~ \frac{|\psi'(r)|^2}{\psi(r)}\leq C
$$ where $C$ is an absolute constant.

Let $$g(X, t)=\psi(\frac{|X|^2}{R^2}).$$  Using the fact that
$|\nabla X|^2=2$,  a straightforward computation shows that,
\begin{eqnarray}\label{e5}
(\Delta-\frac{\partial}{\partial t})g&=&4\psi''\frac{\langle X,
\nabla X\rangle^2}{R^4}+2\psi'\frac{\langle\nabla X, \nabla
X\rangle}{R^2}\nonumber\\ &\geq&-\frac{C_1}{R^2}\nonumber, \\
\frac{|\nabla g|^2}{g}&\leq& \frac{C_2}{R^2}.
\end{eqnarray}

Let $(X(s_i), t(s_i))$ be the point where $g\cdot f$ achieves its
maximum in $\overline{B_{R}(0)}\times [s_i, 0]$. If
$\Sigma_{s_i}\cap B_R(0)=\emptyset$ as $i\to\infty$, then $g\cdot
f\to 0$ as $i\to\infty$. If $\Sigma_{s_i}\cap B_R(0)\neq\emptyset$
as $i\to\infty$, by (\ref{e7}), we know that $f(X, s_i)$ is close
to $\frac{1}{\cos^2\theta(x, s_i)}$ as $i$ large enough, therefore
$f(X, s_i)<e^{h^2}$ for $i$ sufficiently large, since we are
assuming $\delta^2>e^{-h^2}$. We choose $p$ such that $p$ is
sufficiently close to $1$ and keep the condition $1-p>0$. Thus
$f(X, s_i)\leq e^{ph^2}$ as $i\to\infty$. This implies that the
maximum of $g\cdot f$ can not be achieved at $s_i$ as
$i\to\infty$. We can assume that $g\cdot f(X(s_i), t(s_i))>0$. By
the maximum principle, at $(X(s_i), t(s_i))$ we have,
\begin{eqnarray}\label{e1}
&&\nabla(g\cdot f)=0 \nonumber\\ &&\frac{\partial}{\partial
t}(g\cdot f)\geq 0
\end{eqnarray}
and  $$\Delta(g\cdot f)\leq 0. $$ Hence
\begin{eqnarray}\label{e2}(\Delta-\frac{\partial}{\partial
t})g\cdot f\leq 0,\end{eqnarray}
\begin{eqnarray}\label{e11}
\nabla g=-\frac{g}{f}\nabla f.
\end{eqnarray}

Substituting (\ref{e4}) and (\ref{e5}) into (\ref{e2}) and using
(\ref{e11}) twice we get, \allowdisplaybreaks
\begin{eqnarray}\label{e6}
0 &\geq& (\Delta-\frac{\partial}{\partial t}) g\cdot
f=f(\Delta-\frac{\partial}{\partial
t})g+g(\Delta-\frac{\partial}{\partial t})f+2\nabla g\cdot\nabla f
\nonumber
\\ &\geq& -\frac{C_1}{R^2}f-2\frac{|\nabla g|^2}{g} f
+g(\Delta-\frac{\partial}{\partial t})f \nonumber\\
&\geq&-\frac{C_1+2C_2}{R^2}f+ 2g\cdot f|H|^2(1-p)\nonumber\\
&&+g\cdot f(2p|\nabla H|^2+4p^2|H|^2|\nabla
|H||^2-2\frac{|\nabla\cos\theta|^2}{\cos^2\theta})\nonumber\\
&&+2 g\cos^2\theta\nabla f\cdot\nabla\frac{1}{\cos^2\theta} \nonumber\\
&\geq&-\frac{C_1+2C_2}{R^2}f+ 2g\cdot f|H|^2(1-p)\nonumber\\
&&+g\cdot f(2p|\nabla H|^2+4p^2|H|^2|\nabla
|H||^2-2\frac{|\nabla\cos\theta|^2}{\cos^2\theta})\nonumber\\&&-2\cos^2\theta
f\nabla\frac{1}{\cos^2\theta}\cdot\nabla g.
\end{eqnarray} Using the equation (\ref{e11}),
$$\nabla g=g(2\frac{\nabla \cos\theta}{\cos\theta}-p\nabla |H|^2).
$$ Thus,
\begin{eqnarray*}
4gp^2|\nabla |H||^2 |H|^2&=&\frac{|\nabla
g|^2}{g}+4g\frac{|\nabla\cos\theta|^2}{\cos^2\theta}-4\nabla
g\cdot\frac{\nabla\cos\theta}{\cos\theta}.
\end{eqnarray*} Putting this equation into (\ref{e6}),
we get,
\begin{eqnarray*}
0 &\geq& -\frac{C_1+2C_2}{R^2}f+2gf(1-p)|H|^2+2pgf |\nabla
H|^2+\frac{f}{g}|\nabla
g|^2+2gf\frac{|\nabla\cos\theta|^2}{\cos^2\theta} \\ &\geq&
-\frac{C_3}{R^2}f+2gf(1-p)|H|^2.
\end{eqnarray*} This implies that
\begin{eqnarray*}
\frac{C_3}{R^2} &\geq& 2g(1-p)|H|^2 =2gf(1-p)\frac{\cos^2\theta
|H|^2}{e^{p|H|^2}} \\ &\geq& 2gf\delta^2 e^{-ph^2}(1-p)|H|^2.
\end{eqnarray*} By the assumption that $\sup_{t\in (-\infty,
0]}\sup_{\Sigma_t}|A|^2=1$, we have $h^2\leq 2$, so
\begin{eqnarray*}
\frac{C_4}{R^2}\geq \delta^2 2gf(1-p)|H|^2.
\end{eqnarray*}
Since $1-p>0$, we get that,
$$ |H|^2(X(s_i), t(s_i))(g\cdot f)(X(s_i), t(s_i))\leq
\frac{C_4}{(1-p)R^2}.
$$ So,
$$ |H|^2(X(s_i), t(s_i)) f(0, 0)\leq |H|^2(X(s_i), t(s_i))(g\cdot f)(X(s_i), t(s_i))\leq
\frac{C_4}{(1-p)R^2}.
$$ Notice that $f(0, 0)\neq 0$, thus,
$$ |H|^2(X(s_i), t(s_i))\leq
\frac{C_5}{R^2}.
$$  Therefore,
\begin{eqnarray*}
\sup_{B_{\frac{R}{2}}\cap [s_i, 0]} f(X, t)\leq
\frac{1}{\delta^2}e^{p |H|^2(x(s_i), t(s_i))}\leq
\frac{1}{\delta^2} e^{\frac{pC_5}{R^2}}.
\end{eqnarray*} Let $i\to\infty$ then $R\to\infty$ we get that
$$\frac{1}{\delta^2}\geq \sup f\geq e^{ph^2},$$ which contradicts our
assumption because $p$ can be chosen so that it is close to $1$.
This completes the proof of Theorem 2. \hfill Q. E. D.

Now we turn to the the proof of Main Theorem 1. Recall the
evolution equation of the K\"ahler angle in ${\mathbb{C}}^2$ (see
\cite{CL1}), \begin{equation}\label{e3.0}
(\frac{\partial}{\partial t}-\Delta)\cos\alpha=|\overline\nabla
J_{\Sigma_t}|^2\cos\alpha, \end{equation} where $J_{\Sigma_t}$ is
an almost complex structure in a tubular neighborhood of
$\Sigma_t$ in ${\mathbb{C}}^2$ with
\begin{equation}\label{eq2}
\left\{\begin{array}{clcr} J_{\Sigma_t}e_1&=&e_2\\
J_{\Sigma_t}e_2&=&-e_1\\ J_{\Sigma_t}v_1&=&v_2\\
J_{\Sigma_t}v_2&=&-v_1.
\end{array}\right.
\end{equation}

It is showed in \cite{CT} and \cite{CL1} that,
\begin{eqnarray}\label{e3.1}
|\overline\nabla J_{\Sigma_t}|^2 \geq
\frac{1}{2}|H|^2,\end{eqnarray} which implies that
$$ (\frac{\partial}{\partial t}-\Delta)\cos\alpha\geq \frac{1}{2}|H|^2\cos\alpha.
$$
Using the equation (\ref{e3.0}) we can prove one monotonicity
formula along the symplectic mean curvature flow in
${\mathbb{R}}^4$ by the same argument as the one used in the proof
of Proposition 2.1.
\begin{proposition}\label{mono2}
Along the symplectic mean curvature flow $\Sigma_t$ in
${\mathbb{C}}^2$, we have,
\begin{eqnarray*}
&&\frac{\partial}{\partial
t}\left(\int_{\Sigma_t}\frac{1}{\cos\alpha} \rho(F,
t)d\mu_t\right)\nonumber\\&&= -\left(
\int_{\Sigma_t}\frac{1}{\cos\alpha}\rho (F, t)
\left|H+\frac{(F-X_0)^{\perp}}{2(t_0-t)}\right|^2d\mu_t \right.\nonumber \\
&&\left.+ \int_{\Sigma_t}\frac{1}{\cos\alpha}\rho (F,
t)|\overline\nabla J_{\Sigma_t}|^2d\mu_t
+\int_{\Sigma_t}\frac{2}{\cos^3\alpha}\left|\nabla
\cos\alpha\right|^2\rho (F, t)d\mu_t\right).
\end{eqnarray*}
\end{proposition}
By this monotonicity formula we can find a sequence $\{s_i\}$ such
that $s_i\to-\infty$ and
$$\int_{\Sigma_{s_i}\cap B_R(0)}|\overline\nabla
J_{\Sigma_t}|^2\to 0~~~~~~~~~{\rm as}~~~~~~i\to \infty.
$$ By (\ref{e3.1}) we get that,
\begin{eqnarray}\label{e3.2}\lim_{i\to \infty}\max_{\Sigma_{s_i}\cap B_R(0)}|H|^2=0.
\end{eqnarray} We still argue it by contradiction. We assume $\delta>e^{-\frac{h^2}{4}}$
and construct the  function $f=\frac{e^{p|H|^2}}{\cos^2\alpha}$
where $t\in [s_i, 0]$. Due to the inequality (\ref{e3.1}), here
$p$ should be chosen so that $p$ is sufficiently close to $1/2$
and keeps the condition $1/2-p>0$. Using the equation
\begin{eqnarray*}
(\Delta-\frac{\partial}{\partial t})\frac{1}{\cos^2\alpha}
&=&6\frac{|\nabla\cos\alpha|^2}{\cos^4\alpha}+2\frac{|\overline\nabla
J_{\Sigma_t}|^2}{\cos^2\alpha}\\ &\geq& 6
\frac{|\nabla\cos\alpha|^2}{\cos^4\alpha}+\frac{|H|^2}{\cos^2\alpha},
\end{eqnarray*} we obtain that,
\begin{eqnarray}\label{e3.2}
(\Delta-\frac{\partial}{\partial t})f &\geq& f(4p^2|H|^2|\nabla
|H|^2|+2p|\nabla
H|^2+2(1/2-p)|H|^2-2\frac{|\nabla\cos\theta|^2}{\cos^2\theta})
\nonumber\\ &&+2\cos^2\theta\nabla
f\cdot\nabla\frac{1}{\cos^2\theta}.
\end{eqnarray} Similarly we can get,
$$\frac{1}{\delta^2}\geq \sup f\geq e^{ph^2},$$ which
contradicts our assumption that $\delta>e^{-\frac{h^2}{4}}$
because $p$ is close to $1/2$. We leave the details to the reader.

\end{document}